\documentclass[12pt,a4paper]{article}
\usepackage[cp1251]{inputenc} 
\usepackage{amsfonts}

\newcommand{\di}{\displaystyle}
\newcommand{\B}{$\hfill\Box$}
\newcommand{\al}{\alpha}

\newcommand{\ga}{\gamma}
\newcommand{\de}{\delta}
\newcommand{\la}{\lambda}
\newcommand{\om}{\omega}
\newcommand{\ee}{\varepsilon}
\newcommand{\vv}{\varphi}
\newcommand{\iy}{\infty}

\begin{document}

\begin{center}
{\large\bf
Recovering Variable Order Differential Operators with Regular Singularities on Graphs.}\\[0.2cm]
{\bf V.\,Yurko} \\[0.2cm]
\end{center}

\thispagestyle{empty}

{\bf Abstract.} We study inverse spectral problems for ordinary
differential equations with regular singularities on compact
star-type graphs when differential equations have different orders
on diferent edges. As the main spectral characteristics we
introduce and study the so-called Weyl-type matrices which are
generalizations of the Weyl function for the classical
Sturm-Liouville operator. We provide a procedure for constructing
the solution of the inverse problem and prove its uniqueness.



{\bf 1. Introduction. } We study inverse spectral problems for
variable order differential equations with regular singularities
on compact star-type graphs. More precisely, differential
equations have different orders on diferent edges. Boundary value
problems on graphs (spatial networks, trees) often appear in
natural sciences and engineering (see [1-4]). Differential
equations of variable orders on graphs arise in various problems
in mathematics as well as in applications. In particular, we
mention transverse oscillation problems for such structures as
cable-stayed bridges, masts with cable supports and others.

Inverse spectral problems consist in recovering operators from
their spectral characteristics. We pay attention to the most
important nonlinear inverse problems of recovering coefficients of
differential equations (potentials) provided that the structure of
the graph is known a priori.

For {\it second-order} differential operators on compact graphs
inverse spectral problems have been studied fairly completely in
[5-10] and other works. Inverse problems for {\it higher-order}
differential operators on graphs were investigated in [11-12]. We
note that inverse spectral problems for second-order and for
higher-order ordinary differential operators on {\it an interval}
have been studied by many authors (see the monographs [13-18] and
the references therein). Arbitrary order differential operators on
an interval with regular singularities were considered in [19-22].
Variable order differential operators without singularities were
investigated in [23-24]. Variable order differential operators on
graphs with regular singularities have not been studied yet.

In this paper we study the inverse spectral problem for variable
order differential operators with regular singularities on compact
star-type graphs. As the main spectral characteristics in this
paper we introduce and study the so-called Weyl-type matrices
which are generalizations of the Weyl function (m-function) for
the classical Sturm-Liouville operator (see [25]), of the Weyl
matrix for higher-order differential operators on an interval
introduced in [17-18], and generalizations of the Weyl-type
matrices for higher-order differential operators on graphs (see
[11-12]). We show that the specification of the Weyl-type matrices
uniquely determines the coefficients of the differential equation
on the graph, and we provide a constructive procedure for the
solution of the inverse problem from the given Weyl-type matrices.
For studying this inverse problem we develope  the method of
spectral mappings [17-18]. We also essentially use ideas from [19]
on differential equations with regular singularities. The obtained
results are natural generalizations of the well-known results on
inverse problems for differential operators on an interval and on graphs.\\

\medskip
{\bf 2. Weyl-type matrices. } Consider a compact star-type graph $T$ in ${\bf R^\om}$ with the
set of vertices $V=\{v_0,\ldots, v_p\}$ and the set of edges ${\cal E}=\{e_1,\ldots, e_p\},$
where$v_1,\ldots, v_{p}$ are the boundary vertices, $v_0$ is the internal vertex, and
$e_j=[v_{j},v_0],$ $e_1\cap\ldots\cap e_p=\{v_0\}$. Let $l_j$ be the length of
the edge $e_j$. Each edge $e_j\in {\cal E}$ is parameterized by the parameter
$x_j\in [0,l_j]$ such that $x_j=0$ corresponds to the boundary vertices
$v_1,\ldots, v_{p}$, and $x_j=l_j$ corresponds to the internal vertex $v_0$.
An integrable function $Y$ on $T$ may be represented as $Y=\{y_j\}_{j=\overline{1,p}}$,
where the function $y_j(x_j)$ is defined on the edge $e_j$.

Let $n_j$, $j=\overline{1,p},$ be positive integers such that
$n_1\ge n_2\ge\ldots\ge n_p\ge 2.$ Consider the differential
equations on $T$:
$$
y_j^{(n_j)}(x_j)+\di\sum_{\mu=0}^{n_j-2} \Big(\frac{\nu_{\mu
j}}{x_j^{n_j-\mu}}+ q_{\mu j}(x_j)\Big)y_j^{(\mu)}(x_j) =\la
y_j(x_j),\quad x_j\in (0, l_j),\quad j=\overline{1,p},               \eqno(1)
$$
where $\la$ is the spectral parameter, $q_{\mu j}(x_j)$ are
complex-valued integrable functions. We call $q_j=\{q_{\mu
j}\}_{\mu=\overline{0,n_j-2}}$ the potential on the edge $e_j$,
and we call $q=\{q_{j}\}_{j=\overline{1,p}}$ the potential on
the graph $T.$ Let $\{\xi_{kj}\}_{k=\overline{1,n_j}}$ be the
roots of the characteristic polynomial
$$
\de_j(\xi)=\sum_{\mu=0}^{n_j} \nu_{\mu j} \prod_{k=0}^{\mu-1}
(\xi-k), \quad \nu_{n_j,j}:=1,\; \nu_{n_j-1,j}:=0.
$$
For definiteness, we assume that $\xi_{kj}-\xi_{mj}\ne sn_j,
s\in{\bf Z},$ $Re\,\xi_{1j}<\ldots <Re\,\xi_{n_j,j},$
$\xi_{kj}\ne\overline{0,n_j-3}$ (other cases require minor
modifications). We set $\theta_j:=n_j-1-Re\,(\xi_{n_j,j}-\xi_{1j}),$
and assume that the functions $q_{\mu j}^{(\nu)}(x_j),$
$\nu=\overline{0,\mu-1},$ are absolutely continuous, and
$q_{\mu j}^{(\mu)}(x_j)x_j^{\theta_j}\in L(0,l_j).$

\medskip
Fix $j=\overline{1,p}.$ Let $\la=\rho_j^{n_j},$
$\ee_{kj}=\exp(2\pi ik/n_j),$ $k=\overline{0,n_j-1}.$ It is known
that the $\rho_j$ -- plane can be partitioned into sectors
$S_{j\xi}$ of angle $\frac{\pi}{n_j}$
$\Big(\arg\rho_j\in\Big(\frac{\xi\pi}{n_j},\frac{(\xi+1)\pi}{n_j}\Big),$
$\xi=\overline{-n_j,n_j-1}\Big)$ in which the roots $R_{j1},
R_{j2},\ldots, R_{j,n_j}$ of the equation $R^{n_j}-1=0$ can be
numbered in such a way that
$$
Re(\rho_j R_{j1})<Re(\rho_j R_{j2})<\ldots< Re(\rho_j
R_{j,n_j}),\quad\rho_j\in S_{j\xi}.                                     \eqno(2)
$$
Clearly, $R_{jk}=\ee_{\eta_{jk}}$, where $\eta_{j1},\ldots,
\eta_{j,n_j}$ is a permutation of the numbers $0,1,\ldots, n_j-1,$
depending on the sector. Let us agree that
$$
\rho_j^{\mu}=\exp(\mu(\ln|\rho_j|+i\arg\,\rho_j)),\;
\arg\,\rho_j\in(-\pi,\pi], \quad R_{jk}^{\mu}= \exp(2\pi
i\mu\eta_{jk}/n_j).
$$
Let the numbers $c_{kj0},\; k=\overline{1,n_j},$ be such that
$$
\prod_{k=1}^{n_j} c_{kj0}=
\Big(\det[\xi_{kj}^{\nu-1}]_{k,\nu=\overline{1,n_j}}\Big)^{-1}.
$$
Then the functions
$$
C_{kj}(x_j,\la)=x_j^{\xi_{kj}} \sum_{\mu=0}^\iy c_{kj\mu} (\rho_j
x_j)^{n_j\mu}, \quad c_{kj\mu}=c_{kj0}\Big(\prod_{s=1}^\mu
\de_j(\xi_{kj}+sn_j)\Big)^{-1},
$$
are solutions of the differential equation in the case when
$q_{\mu j}(x_j)\equiv 0,$ $\mu=\overline{0,n_j-2}.$ Moreover,
$\det[C_{kj}^{(\nu-1)}(x_j,\la)]_{k,\nu=\overline{1,n_j}}\equiv
1.$ Denote $\rho^*=\max\limits_{j=\overline{1,p}}
\Big(2n_j\max\limits_{\mu=\overline{0,n_j-2}} \|q_{\mu
j}\|_{L(0,l_j)}\Big).$ In [19] we constructed
 special fundamental systems of solutions
$\{S_{kj}(x_j,\la)\}_{k=\overline{1,n_j}}$ and
$\{E_{kj}(x_j,\rho_j)\}_{k=\overline{1,n_j}}$ of equation (1) on
the edge $e_j$, possessing the following properties.

1) For each $x_j\in(0,l_j],$ the functions
$S_{kj}^{(\nu)}(x_j,\la),\;\nu=\overline{0,n_j-1},$ are entire in
$\la.$ For each fixed $\la,$ and $x_j\to 0,$
$$
S_{kj}(x_j,\la)\sim c_{kj0} x_j^{\xi_{kj}},\quad
(S_{kj}(x_j,\la)-C_{kj}(x_j,\la))x_j^{-\xi_{kj}}=o(x_j^{\xi_{n_j,j}-\xi_{1j}}).
$$
Moreover,
$\det[S_{kj}^{(\nu-1)}(x_j,\la)]_{k,\nu=\overline{1,n_j}}\equiv
1,$ and $|S_{kj}^{(\nu)}(x_j,\la)|\le C|x_j^{\xi_{kj}-\nu}|,\;
|\rho_j|x_j\le 1.$ Here and below, we shall denote by the same
symbol $C$ various positive constants in the estimates independent
of $\la$ and $x_j$.

2) For each $x_j>0$ and for each sector $S_{j\xi}$ with property
(2), the functions $E_{kj}^{(\nu)}(x_j,\rho_j),$
$\nu=\overline{0,n_j-1},$ are regular with respect to $\rho_j\in
S_{j\xi},\; |\rho_j|>\rho^*$, and continuous for $\rho_j\in
\overline{S_{j\xi}},\; |\rho_j|\ge \rho^*$. Moreover,
$$
|E_{kj}^{(\nu)}(x_j,\rho_j)(\rho_j R_{jk})^{-\nu}\exp(-\rho_j
R_{jk} x_j)-1| \le C(|\rho_j|x_j),\quad \rho_j\in
\overline{S_{j\xi}},\quad |\rho_j|x_j\ge 1.
$$

3) The relation
$$
E_{kj}(x_j,\rho_j)=\sum_{\mu=1}^{n_j} b_{kj\mu}(\rho_j)S_{\mu j}(x_j,\la), \eqno(3)
$$
holds, where
$$
b_{kj\mu}(\rho_j)= b_{\mu j}^0 R_{jk}^{\xi_{\mu
j}}\rho_j^{\xi_{\mu j}}[1],\quad b_{\mu j}^0\ne 0,\quad
\rho_j\in\overline{S_{j\xi}},\quad \rho_j\to\iy,                        \eqno(4)
$$
$$
\prod_{\mu=1}^{n_j} b_{\mu
j}^0=\det[R_{jk}^{\nu-1}]_{k,\nu=\overline{1,n_j}}
\Big(\det[R_{jk}^{\xi_{\mu
j}}]_{k,\mu=\overline{1,n_j}}\Big)^{-1},
$$
where $[1]=1+O(\rho^{-1}).$

Note that the asymptotical formula (4) is the most important and
nontrivial property of these solutions. This property allows one
to study both direct and inverse problems for arbitrary order
differential operators with regular singularities (see [20-22]).

Consider the linear forms
$$
U_{j\nu}(y_j)=\sum_{\mu=0}^{\nu}\ga_{j\nu\mu}y_j^{(\mu)}(l_j), \;
j=\overline{1,p},\; \nu=\overline{0,n_j-1},
$$
where $\ga_{j\nu\mu}$ are complex numbers,
$\ga_{j\nu}:=\ga_{j\nu\nu}\ne 0.$ The linear forms $U_{j\nu}$ will
be used in matching conditions at the internal vertex $v_0$ for
boundary value problems and for the correspondung special
solutions of equation (1).

\medskip
Denote $\langle n\rangle:=(|n|+n)/2,$ i.e. $\langle n\rangle=n$
for $n\ge 0$, and $\langle n\rangle=0$ for $n\le 0.$ Fix
$s=\overline{1,p},\; k=\overline{1,n_s-1}.$ Let $\Psi_{sk}=
\{\psi_{skj}\}_{j=\overline{1,p}}$ be solutions of equation (1) on
the graph $T$ under the boundary conditions
$$
\psi_{sks}(x_s,\la)\sim c_{ks0}x_s^{\xi_{ks}},\quad x_s\to 0,     \eqno(5)
$$
$$
\psi_{skj}(x_j,\la)=O(x_j^{\xi_{\langle n_j-k-1\rangle+2,j}}),
\quad x_j\to 0,\quad j=\overline{1,p},\;j\ne s,                  \eqno(6)
$$
and the matching conditions at the vertex $v_0$:
$$
U_{1\nu}(\psi_{sk1})=U_{j\nu}(\psi_{skj}), \;\;
j=\overline{2,p},\; \nu=\overline{0,k-1},\;n_j>\nu+1,            \eqno(7)
$$
$$
\di\sum_{j=1,\,n_j>\nu}^{p} U_{j\nu}(\psi_{skj})=0,
\quad \nu=\overline{k,n_s-1}.                                    \eqno(8)
$$
In particular, if $n_j>k,$ then condition (6) takes the form
$$
\psi_{skj}(x_j,\la)=O(x_j^{\xi_{n_j-k+1,j}}), \quad x_j\to 0,
\quad j=\overline{1,p},\;j\ne s,
$$
and if $n_j\le k,$ then condition (6) takes the form
$$
\psi_{skj}(x_j,\la)=O(x_j^{\xi_{2j}}), \quad x_j\to 0,\quad
j=\overline{1,p},\;j\ne s.
$$
Matching conditions (7)-(8) are generalizations of classical
matching conditions for higher-order differential operators on
graphs, and matching conditions for variable order differential
operators on graphs [23-24]. The function $\Psi_{sk}$ is called
the Weyl-type solution of order $k$ with respect to the boundary
vertex $v_s$. Define additionally
$\psi_{sns}(x_s,\la):=S_{ns}(x_s,\la).$

Using the fundamental system of solutions $\{S_{\mu j}(x_j,\la)\}$
on the edge $e_j$, one can write
$$
\psi_{skj}(x_j,\la)=\di\sum_{\mu=1}^{n_j} M_{skj\mu}(\la)
S_{\mu j}(x_j,\la),\quad j=\overline{1,p},\quad k=\overline{1,n_s-1}, \eqno(9)
$$
where the coefficients $M_{skj\mu}(\la)$ do not depend on $x_j.$
It follows from (9) and the boundary condition (5) for the
Weyl-type solutions that
$$
\psi_{sks}(x_s,\la)=S_{ks}(x_s,\la)+\sum_{\mu=k+1}^{n_s}
M_{sk\mu}(\la)S_{\mu s}(x_s,\la), \quad
M_{sk\mu}(\la):=M_{sks\mu}(\la).  \eqno(10)
$$
We introduce the matrices $M_{s}(\la),\;s=\overline{1,p},$ as
follows:
$$
M_{s}(\la)=[M_{sk\mu}(\la)]_{k,\mu=\overline{1,n_s}},\quad
M_{sk\mu}(\la):=\de_{k\mu}\quad \mbox{for}\quad k\ge \nu.
$$
The matrix $M_s(\la)$ is called the Weyl-type matrix with respect
to the boundary vertex $v_s$. The inverse problem is formulated as
follows. Fix $w=\overline{2,p}.$

\smallskip
{\bf Inverse problem 1.} Given $\{M_{s}(\la)\},\;
s=\overline{1,p}\setminus w$, construct $q$ on $T.$

\smallskip
We note that the notion of the Weyl-type matrices $M_s$ is a
generalization of the notion of the Weyl function (m-function) for
the classical Sturm-Liouville operator ([15, 25])  and is a
generalization of the notion of Weyl matrices introduced in [11,
12, 17, 18, 20] for higher-order differential operators on an
interval and on graphs. Thus, Inverse Problem 1 is a
generalization of the well-known inverse problems for differential
operators on an interval and on graphs.

\smallskip
We also note that in Inverse problem 1 we do not need to specify
all matrices $M_s(\lambda),$ $s=\overline{1,p}$; one of them can
be omitted. This last fact was first noticed in [6], where the
inverse problem was solved for the Sturm-Liouville operators on an
arbitrary tree.

\smallskip
In section 3 properties of the Weyl-type solutions  and the
Weyl-type matrices are studied. Section 4 is devoted to the
solution of auxiliary inverse problems of recovering the potential
on a fixed edge. In section 5 we study Inverse Problem 1. For this
inverse problem we provide a constructive procedure for the
solution and prove its uniqueness.\\

{\bf 3. Properties of spectral characteristics.} Fix
$s=\overline{1,p},\; k=\overline{1,n_s-1}.$ Substituting (9) into
boundary and matching conditions (5)-(8), we obtain a linear
algebraic system with respect to $M_{skj\mu}(\la).$ Solving this
system by Cramer's rule one gets
$M_{skj\mu}(\la)=\Delta_{skj\mu}(\la)/\Delta_{sk}(\la),$ where the
functions $\Delta_{skj\mu}(\la)$ and $\Delta_{sk}(\la)$ are entire
in $\la.$ Thus, the functions $M_{skj\mu}(\la)$ are meromorphic in
$\la,$ and consequently, the Weyl-type solutions and the Weyl-type
matrices are meromorphic in $\la.$ In particular,
$$
M_{sk\mu}(\la)= \frac{\Delta_{sk\mu}(\la)}{\Delta_{sk}(\la)},\quad
k\le\mu,
$$
where $\Delta_{sk\mu}(\la):=\Delta_{sks\mu}(\la).$ We note that
the function $\Delta_{sk}(\la)$ is the characteristic function for
the boundary value problem $L_{sk}$ for equation (1) under the
conditions
$$
y_{s}(x_s)=O(x_s^{\xi_{k+1,s}}),\; x_s\to 0,\qquad y_{j}(x_j)=
O(x_j^{\xi_{\langle n_j-k-1\rangle+2,j}}) ,\; x_j\to 0,\;
j=\overline{1,p},\;j\ne s,
$$
$$
U_{1\nu}(y_{1})=U_{j\nu}(y_{j}),\;
j=\overline{2,p},\;\nu=\overline{0,k-1},\;n_j>\nu+1,
$$
$$
\sum_{j=1,\;n_j>\nu}^{p}
U_{j\nu}(y_{j})=0,\;\nu=\overline{k,n_s-1}.
$$
Zeros of $\Delta_{sk}(\la)$ coincide with the eigenvalues of
$L_{sk}.$ Denote
$$
\Omega_{jk}=\det[R_{jl}^{\xi_{\mu j}}]_{l,\mu=\overline{1,k}},\quad
\Omega_{j0}=1, \quad \om_{jk}:=\frac{\Omega_{j,k-1}}{\Omega_{jk}},\;
k=\overline{1,n_j}.
$$

{\bf Lemma 1. }{\it Fix $j=\overline{1,p},$ and fix a sector $S_{j\xi}$
with property (2).

1) Let $k=\overline{1,n_j-1},$ and let $y_j(x_j,\la)$ be a solution
of equation (1) on the edge $e_j$ under the condition
$$
y_j(x_j,\la)=O(x_j^{\xi_{k+1,j}}),\quad x_j\to 0.                      \eqno(10)
$$
Then for $x_j\in(0,l_j],\;\nu=\overline{0,n_j-1},\;
\rho_j\in S_{j\xi},\; |\rho_j|\to\iy,$
$$
y_j^{(\nu)}(x_j,\la)=\sum_{\mu=k+1}^{n_j} A_{\mu j}(\rho_j)
(\rho_j R_{j\mu})^{\nu}\exp(\rho_j R_{j\mu} x_j)[1],                  \eqno(11)
$$
where the coefficients $A_{\mu j}(\rho_j)$ do not depend on $x_j$.
Here and below we assume that $\arg\rho_j=const,$ when $|\rho_j|\to\iy.$

2) Let $k=\overline{1,n_j},$ and let $y_j(x_j,\la)$ be a solution
of equation (1) on the edge $e_j$ under the condition
$$
y_j(x_j,\la)\sim c_{kj0}x_j^{\xi_{kj}},\; x_j\to 0.                   \eqno(12)
$$
Then for $x_j\in(0,l_j],\;\nu=\overline{0,n_j-1},\;
\rho_j\in S_{j\xi},\;|\rho_j|\to\iy,$
$$
y_j^{(\nu)}(x_j,\la)=\frac{\om_{jk}}{\rho_j^{\xi_{kj}}}
(\rho_j R_{jk})^{\nu}\exp(\rho_j R_{jk} x_j)[1]
+\sum_{\mu=k+1}^{n_j} B_{\mu j}(\rho_j)
(\rho_j R_{j\mu})^{\nu}\exp(\rho_j R_{j\mu} x_j)[1],                  \eqno(13)
$$
where the coefficients $B_{\mu j}(\rho_j)$ do not depend on $x_j$.}

\smallskip
{\bf Proof. } It follows from (10) that
$$
y_j(x_j,\la)=\sum_{\mu=k+1}^{n_j} a_{\mu j}(\la)S_{\mu j}(x_j,\la).   \eqno(14)
$$
Using the fundamental system of solutions
$\{E_{kj}(x_j,\rho_j)\}_{k=\overline{1,n_j}}$, one can write
$$
y_j(x_j,\la)=\sum_{m=1}^{n_j} A_{mj}(\rho_j)E_{mj}(x_j,\rho_j).       \eqno(15)
$$
By virtue of (3), we calculate
$$
y_j(x_j,\la)=\sum_{m=1}^{n_j} A_{mj}(\rho_j)\sum_{\mu=1}^{n_j}
b_{mj\mu}(\rho_j) S_{\mu j}(x_j,\la)=\sum_{\mu=1}^{n_j} S_{\mu j}(x_j,\la)
\sum_{m=1}^{n_j} A_{mj}(\rho_j) b_{mj\mu}(\rho_j).
$$
Comparing this relation with (14), we obtain
$$
\sum_{m=1}^{n_j} A_{mj}(\rho_j) b_{mj\mu}(\rho_j)=0,
\quad \mu=\overline{1,k}.                                            \eqno(16)
$$
We consider (16) as a linear algebraic system with respect to
$A_{j}(\rho_j), A_{2j}(\rho_j),\ldots, A_{kj}(\rho_j).$ Solving
this system by Cramer's rule and taking (4) into account we get
$$
A_{mj}(\rho_j)=\sum_{\mu=k+1}^{n_j}(\al_{m\mu j}
+O(\rho_j^{-1}))A_{\mu j}(\rho_j),\quad m=\overline{1,k},            \eqno(17)
$$
where $\al_{m\mu j}$ are constants. Substituting (17) into (15)
and using (2) we arrive at (11). Relations (13) are proved
analogously by using (12) instead of (10).
\B

\smallskip
Now we are going to study the asymptotic behavior
of the Weyl-type solutions.

\smallskip
{\bf Lemma 2. }{\it Fix $s=\overline{1,p},\;k=\overline{1,n_s},$
and fix a sector $S_{s\xi}$ with property (2). For $x_s\in(0,l_s),\;
\nu=\overline{0,n_s-1},$ the following asymptotic formula holds}
$$
\psi_{sks}^{(\nu)}(x_s,\la)=\frac{\om_{sk}}{\rho_s^{\xi_{ks}}}\,
(\rho_s R_{sk})^{\nu}\exp(\rho_s R_{sk} x_s)[1],
\quad \rho_s\in S_{s\xi},\;|\rho_s|\to\iy.                             \eqno(18)
$$

{\bf Proof. } For $k=n_s$, (18) follows from Lemma 1. Fix
$s=\overline{1,p},\; k=\overline{1,n_s-1}.$ Using Lemma 1 and
boundary conditions for $\Psi_{sk}$ we get the following asymptotic
formulae for $|\la|\to\iy,$ inside the corresponding sectors:
$$
\psi_{sks}^{(\nu)}(x_s,\la)=\frac{\om_{sk}}{\rho_s^{\xi_{ks}}}
\,(\rho_s R_{sk})^{\nu}\exp(\rho_s R_{sk} x_s)[1]+
\sum_{\mu=k+1}^{n_s} A_{\mu s}^{sk}(\rho_s)
(\rho_s R_{s\mu})^{\nu}\exp(\rho_s R_{s\mu} x_s)[1],\; x_s\in(0,l_s],  \eqno(19)
$$
$$
\psi_{skj}^{(\nu)}(x_j,\la)=\sum_{\mu=\langle n_j-k+1\rangle +2}^{n_j}
A_{\mu j}^{sk}(\rho_j)(\rho_j R_{j\mu})^{\nu}\exp(\rho_j R_{j\mu} x_j)[1],
\;j=\overline{1,p}\setminus s,\;  x_j\in(0,l_j].                       \eqno(20)
$$
Substituting (19)-(20) into matching conditions (7)-(8) for
$\Psi_{sk}$, we obtain the linear algebraic system with respect to
$A_{\mu j}^{sk}(\rho_j).$ Solving this system by Cramer's rule,
we obtain in particular,
$$
A_{\mu s}^{sk}(\rho_s)
=O(\rho_s^{-\xi_{ks}}\exp(\rho_s(R_{sk}-R_{s\mu})l_s)).              \eqno(21)
$$
Substituting (21) into (19) we arrive at (18).
\B

\smallskip
It follows from the proof of Lemma 2 that one can also get the
asymptotics for $\psi_{skj}^{(\nu)}(x_j,\la),$ $j\ne s$; but for
our purposes only (18) is needed.\\

{\bf 4. Auxiliary inverse problems. } In this section we consider
auxiliary inverse problems of recovering differential operator on
each fixed edge. Fix $s=\overline{1,p},$ and consider the
following auxiliary inverse problem on the edge $e_s$.

\smallskip
{\bf IP(s). } Given the matrix $M_s$,
construct the potential $q_{s}$ on the edge $e_s$.

\smallskip
In this inverse problem we construct the potential only
on the edge $e_s$, but the Weyl-type matrix $M_s$ brings a global
information from the whole graph. In other words, this problem
is not a local inverse problem related only to the edge $e_s$.

\smallskip
{\bf Theorem 1. }{\it Fix $s=\overline{1,p}.$ The specification of the
Weyl-type matrix $M_s$ uniquely determines the potential $q_s$ on the
edge $e_s$.}

\smallskip
We omit the proof since it is similar to that in [18, Ch.2].
Moreover, using the method of spectral mappings and the
asymptotics (18) for the Weyl-type solutions, one can get a
constructive procedure for the solution of the inverse problem
$IP(s)$. It can be obtained by the same arguments as for $n$-th
order differential operators on a finite interval (see [18, Ch.2]
for details). Note that like in [18], the nonlinear inverse
problem $IP(s)$ is reduced to the solution of a linear equation in
the corresponding Banach space of sequences.

\smallskip
Fix $j=\overline{1,p}.$ Let $\vv_{jk}(x_j,\la),$ $k=\overline{1,n_j},$
be solutions of equation (1) on the edge $e_j$ under the conditions
$$
\vv_{kj}^{(\nu-1)}(l_j,\la)=\de_{k\nu},\;\nu=\overline{1,k},
\qquad \vv_{kj}(x_j,\la)=O(x_j^{\xi_{n_j-k+1,j}}),\;x_j\to 0.
$$
We introduce the matrix
$m_j(\la)=[m_{jk\nu}(\la)]_{k,\nu=\overline{1,n_j}},$
where $m_{jk\nu}(\la):=\vv^{(\nu-1)}_{jk}(l_j,\la).$
The matrix $m_j(\la)$ is called the Weyl-type matrix with
respect to the internal vertex $v_0$ and the edge $e_j$.

\smallskip
{\bf IP[j]. } Given the matrix $m_j$, construct $q_{j}$ on the edge $e_j$.

This inverse problem is the classical one, since it is the inverse
problem of recovering a higher-order differential equation on a
finite interval from its Weyl-type matrix. This inverse problem
has been solved in [18], where the uniqueness theorem for this
inverse problem is proved. Moreover, in [18] an algorithm for the
solution of the inverse problem $IP[j]$ is given, and necessary
and sufficient conditions for the solvability of this inverse
problem are provided.\\

{\bf 5. Solution of Inverse Problem 1. } In this section we obtain
a constructive procedure for the solution of Inverse problem 1 and
prove its uniqueness. First we prove an auxiliary assertion.

\smallskip
{\bf Lemma 3. }{\it Fix $j=\overline{1,p}.$
Then for each fixed $s=\overline{1,p}\setminus j,$}
$$
m_{j1\nu}(\la)=\di\frac{\psi_{s1j}^{(\nu-1)}(l_j,\la)}{\psi_{s1j}(l_j,\la)},
\quad \nu=\overline{2,n_j},                                                \eqno(22)
$$
$$
m_{jk\nu}(\la)=\di\frac{\det[\psi_{s\mu j}(l_j,\la),\ldots,
\psi_{s\mu j}^{(k-2)}(l_j,\la),
\psi_{s\mu j}^{(\nu-1)}(l_j,\la)]_{\mu=\overline{1,k}}}
{\det[\psi_{s\mu j}^{(\xi-1)}(l_j,\la)]_{\xi,\mu=\overline{1,k}}}\,,
\;2\le k<\nu\le n_j.                                                      \eqno(23)
$$

{\bf Proof. } Denote
$$
w_{js}(x_j,\la):=\di\frac{\psi_{s1j}(x_j,\la)}{\psi_{s1j}(l_j,\la)}\,.
$$
The function $w_{js}(x_j,\la)$ is a solution of equation (1) on the edge
$e_j$, and $w_{js}(l_j,\la)=1.$ Moreover, by virtue of the boundary
conditions on $\Psi_{s1}$, one has $w_{js}(x_j,\la)=O(x_j^{\xi_{n_j,j}}),$
$x_j\to 0.$ Hence, $w_{js}(x_j,\la)\equiv\vv_{1j}(x_j,\la),$
i.e.
$$
\vv_{1j}(x_j,\la)=\di\frac{\psi_{s1j}(x_j,\la)}{\psi_{s1j}(l_j,\la)}.     \eqno(24)
$$
Similarly, we calculate
$$
\vv_{kj}(x_j,\la)=\di\frac{\det[\psi_{s\mu j}(l_j,\la),\ldots,
\psi_{s\mu j}^{(k-2)}(l_j,\la),\psi_{s\mu j}(x_j,\la)]_{\mu=\overline{1,k}}}
{\det[\psi_{s\mu j}^{(\xi-1)}(l_j,\la)]_{\xi,\mu=\overline{1,k}}}\,,
\quad k=\overline{2,n_j-1}.                                               \eqno(25)
$$
Since $m_{jk\nu}(\la)=\vv_{kj}^{(\nu-1)}(l_j,\la),$ it follows
from (24) that (22) holds. Similarly, (23) follows from (25).
\B

\medskip
Now we are going to obtain a constructive procedure for the
solution of Inverse problem 1. Our plan is the following.

{\it Step 1. } Let the Weyl-type matrices $\{M_{s}(\la)\},\;
s=\overline{1,p}\setminus w$, be given. Solving the inverse problem
$IP(s)$ for each fixed $s=\overline{1,p}\setminus w,$ we find the
potentials $q_{s}$ on the edges $e_s$, $s=\overline{1,p}\setminus w$.

{\it Step 2. } Using the knowledge of the potential on the edges
$e_s$, $s=\overline{1,p}\setminus w$, we construct the Weyl-type
matrix $m_{w}(\la)$.

{\it Step 3. } Solving the inverse problem $IP[w]$, we find
the potential $q_{w}$ on $e_{w}$.

\smallskip
Steps 1 and 3 have been already studied in Section 3. It remains to
fulfil Step 2. For this purpose it is convenient to divide differential
equations into $m$ groups with equal orders. More precisely, let
$\om_1>\om_2>\ldots >\om_m>\om_{m+1}=1,$ $n_{p_{j-1}+1}=\ldots=n_{p_j}
:=\om_j,$ $j=\overline{1,m},$ $0=p_0<p_1<\ldots <p_m=p.$
Take $N$ such that $p_N=w.$

Suppose that Step 1 is already made, and we found the potentials
$q_{s}$, $s=\overline{1,p}\setminus p_N$, on the edges $e_s$,
$s=\overline{1,p}\setminus p_N$. Then we calculate the functions
$S_{kj}(x_j,\la),$ $j=\overline{1,p}\setminus p_N;$
here $k=\overline{1,\om_i}$ for $j=\overline{p_{i-1}+1,p_i}.$

Fix $s=\overline{1,p_1}$ (if $N>1$), and $s=\overline{1,p_1-1}$ (if $N=1$).
All calculations below will be made for this fixed $s.$
Our goal now is to construct the Weyl-type matrix $m_{p_N}(\la).$
According to (22)-(23), in order to construct $m_{p_N}(\la)$ we have
to calculate the functions
$$
\psi_{skp_N}^{(\nu)}(l_{p_N},\la),
\quad k=\overline{1,\om_N-1},\; \nu=\overline{0,\om_N-1}.                        \eqno(26)
$$
We will find the functions (26) by the following steps.

\smallskip
1) Using (10) we construct the functions
$$
\psi_{sks}^{(\nu)}(l_s,\la),\;k=\overline{1,\om_N-1},\;\nu=\overline{0,\om_1-1}, \eqno(27)
$$
by the formula
$$
\psi_{sks}^{(\nu)}(l_s,\la)=S_{ks}^{(\nu)}(l_s,\la)+
\di\sum_{\mu=k+1}^{\om_1} M_{sk\mu}(\la)S_{\mu s}^{(\nu)}(l_s,\la).              \eqno(28)
$$

2) Consider a part of the matching conditions (7) on $\Psi_{sk}$. More
precisely, let $\xi=\overline{N,m},\;k=\overline{\om_{\xi+1},\om_\xi-1},\;
l=\overline{\xi,m},\;j=\overline{1,p_l-1}.$ Then, in particular, (7) yields
$$
U_{p_l,\nu}(\psi_{skp_l})=U_{j\nu}(\psi_{skj}),
\quad \nu=\overline{\om_{l+1}-1,\min(k-1,\om_l-2)}.                              \eqno(29)
$$
Since the functions (27) are known, it follows from (29) that
one can calculate the functions
$$
\psi_{skj}^{(\nu)}(l_j,\la),\;\xi=\overline{N,m},\;
k=\overline{\om_{\xi+1},\om_\xi-1},\; l=\overline{\xi,m},\;
j=\overline{1,p_l},\; \nu=\overline{\om_{l+1}-1,\min(k-1,\om_l-2)}.             \eqno(30)
$$
In particular we found the functions (26) for $\nu=\overline{0,k-1}.$

\smallskip
3) It follows from (9) and the boundary conditions on $\Psi_{sk}$ that
$$
\psi_{skj}^{(\nu)}(l_j,\la)=\di\sum_{\mu=\max(\om_l-k+1)}^{\om_l}
M_{skj\mu}(\la)C_{\mu j}^{(\nu)}(l_j,\la),                                     \eqno(31)
$$
$$
k=\overline{1,\om_1-1},\;l=\overline{1,m},\;
j=\overline{p_{l-1}+1,p_l}\setminus s,\;\nu=\overline{0,\om_l-1}.
$$
We consider only a part of relations (31). More precisely, let
$\xi=\overline{N,m},\; k=\overline{\om_{\xi+1},\om_\xi-1},\;
l=\overline{1,m},\;j=\overline{p_{l-1}+1,p_l},\; j\ne p_N,\;
j\ne s,\; \nu=\overline{0,\min(k-1,\om_l-2)}.$ Then
$$
\di\sum_{\mu=\max(\om_l-k+1)}^{\om_l}M_{skj\mu}(\la)C_{\mu j}^{(\nu)}(l_j,\la)
=\psi_{skj}^{(\nu)}(l_j,\la),\quad \nu=\overline{0,\min(k-1,\om_l-2)}.        \eqno(32)
$$
For this choice of parameters, the right-hand side in (32) are known,
since the functions (30) are known. Relations (32) form a linear algebraic
system $\sigma_{skj}$ with respect to the coefficients $M_{skj\mu}(\la).$
Solving the system $\sigma_{skj}$ by Cramer's rule we find the functions
$M_{skj\mu}(\la).$  Substituting them into (31), we calculate the functions
$$
\psi_{skj}^{(\nu)}(l_j,\la),\quad k=\overline{1,\om_N-1},\; l=\overline{1,m},
\;j=\overline{p_{l-1}+1,p_l}\setminus p_N,\; \nu=\overline{0,\om_l-1}.        \eqno(33)
$$
Note that for $j=s$ these functions were found earlier.

\smallskip
4) Let us now use the generalized Kirchhoff'sconditions (8) for $\Psi_{sk}$.
Since the functions (33) are known, one can construct by (8) the functions
(26) for $k=\overline{1,\om_N-1},\;\nu=\overline{k,\om_N-1}.$ Thus, the
functions (26) are known for $k=\overline{1,\om_N-1},\; \nu=\overline{0,\om_N-1}.$

\smallskip
Since the functions (26) are known, we construct the Weyl-type matrix
$m_{p_N}(\la)$ via (22)-(23) for $j=p_N.$ Thus, we have obtained the solution of
Inverse problem 1 and proved its uniqueness, i.e. the following assertion holds.

\medskip
{\bf Theorem 2. }{\it The specification of the Weyl-type matrices $M_s(\la),$
$s=\overline{1,p}\setminus p_N$, uniquely determines the potential $q$ on $T.$
The solution of Inverse problem 1 can be obtained by the following algorithm.}

{\bf Algorithm 1. }{\it Given the Weyl-type matrices
$M_s(\la),$ $s=\overline{1,p}\setminus p_N$.

1) Find $q_{s}$, $s=\overline{1,p}\setminus p_N$, by solving the
inverse problem $IP(s)$ for each fixed $s=\overline{1,p}\setminus p_N$.

2) Calculate $C_{kj}^{(\nu)}(l_j,\la),\; j=\overline{1,p}\setminus p_N$; here
$k=\overline{1,\om_i},\;\nu=\overline{0,\om_i-1}$ for $j=\overline{p_{i-1}+1,p_i}.$

3) Fix $s=\overline{1,p_1}$ (if $N>1$), and $s=\overline{1,p_1-1}$ (if $N=1$).
All calculations below will be made for this fixed $s.$
Construct the functions (27) via (28).

4) Calculate the functions (30) using (29).

5) Find the functions $M_{skj\mu}(\la),$ by solving the linear algebraic
systems $\sigma_{skj}$.

6) Construct the functions (26) using (8).

7) Calculate the Weyl-type matrix $m_{p_N}(\la)$ via (22)-(23) for $j=p_N$.

8) Construct the potential $q_{p_N}$ on the edge $e_{p_N}$
by solving the inverse problem $IP[j]$ for $j=p_N$.}

\medskip
{\bf Acknowledgment.} This work was supported by Grant 1.1436.2014K of the Russian 
Ministry of Education and Science and by Grant 13-01-00134 of Russian Foundation for 
Basic Research.

\begin{center}
{\bf REFERENCES}
\end{center}
\begin{enumerate}
\item[{[1]}] J. Langese, G. Leugering and J. Schmidt. Modelling, analysis and control
     of dynamic elastic multi-link structures. Birkh\"auser, Boston, 1994.
\item[{[2]}] T. Kottos and U. Smilansky. Quantum chaos on graphs. Phys. Rev. Lett. 79 (1997), 4794-4797.
\item[{[3]}] P. Kuchment. Quantum graphs. Some basic structures. Waves Random Media 14 (2004), S107-S128.
\item[{[4]}] Yu. Pokornyi and A. Borovskikh. Differential equations on  networks (geometric graphs).
     J. Math. Sci. (N.Y.) 119, no.6 (2004), 691-718.
\item[{[5]}] M.I. Belishev. Boundary spectral inverse problem on a class of graphs (trees)
     by the BC method. Inverse Problems 20 (2004), 647-672.
\item[{[6]}] V.A. Yurko. Inverse spectral problems for Sturm-Liouville operators on graphs.
     Inverse Problems 21 (2005), 1075-1086.
\item[{[7]}] B.M. Brown and R. Weikard. A Borg-Levinson theorem for trees.
     Proc. R. Soc. Lond. Ser. A Math. Phys. Eng. Sci. 461, no.2062 (2005), 3231-3243.
\item[{[8]}] V.A. Yurko. Inverse problems for Sturm-Liouville operators on bush-type graphs.
     Inverse Problems 25, no.10 (2009), 105008, 14pp.
\item[{[9]}] V.A. Yurko. An inverse problem for Sturm-Liouville operators on A-graphs.
     Applied Math. Lett. 23, no.8 (2010), 875-879.
\item[{[10]}] V.A. Yurko. Inverse spectral problems for differential operators on arbitrary
     compact graphs. Journal of Inverse and Ill-Posed Proplems 18, no.3 (2010), 245-261.
\item[{[11]}] V.A. Yurko. An inverse problem for higher-order differential operators on
     star-type graphs. Inverse Problems 23, no.3 (2007), 893-903.
\item[{[12]}] V.A. Yurko. Inverse problems for differential of any order on trees.
     Matemat. Zametki 83,no.1 (2008), 139-152; English  transl. in Math. Notes 83, no.1 (2008), 125-137.
\item[{[13]}] V.A. Marchenko. Sturm-Liouville operators and their applications.
     "Naukova Dumka",  Kiev, 1977;  English  transl., Birkh\"auser, 1986.
\item[{[14]}] B.M. Levitan. Inverse Sturm-Liouville problems. Nauka, Moscow, 1984;
     English transl., VNU Sci.Press, Utrecht, 1987.
\item[{[15]}] G. Freiling and V.A. Yurko. Inverse Sturm-Liouville Problems and their Applications.
     NOVA Science Publishers, New York, 2001.
\item[{[16]}] R. Beals, P. Deift and C. Tomei. Direct and Inverse Scattering on the Line,
     Math. Surveys and Monographs, v.28. Amer. Math. Soc. Providence: RI, 1988.
\item[{[17]}] V.A. Yurko. Inverse Spectral Problems for Differential Operators and their Applications.
     Gordon and Breach, Amsterdam, 2000.
\item[{[18]}] V.A. Yurko. Method of Spectral Mappings in the Inverse Problem Theory,
     Inverse and Ill-posed Problems Series. VSP, Utrecht, 2002.
\item[{[19]}] V.A. Yurko. Inverse  problem for differential equations
     with a singularity. Differ. Uravneniya, 28, no.8 (1992), 1355-1362
     (Russian); English transl. in Diff. Equations, 28 (1992),
     1100-1107.
\item[{[20]}] V.A. Yurko. On higher-order differential operators with
     a singularity.  Matem. Sbornik, 186, no.6 (1995), 133-160 (Russian);
     English transl. in Sbornik; Mathematics 186, no.6 (1995), 901-928.
\item[{[21]}] V.A. Yurko. Inverse spectral problems for higher-order
     differential operators with a singularity. Journal of Inverse and
     Ill-Posed Problems, 10, no.4 (2002), 413-425.
\item[{[22]}] V.A. Yurko. Higher-order differential equations having
     a singularity in an interior point. Results in Mathematics, 42, no.1-2 (2002),
     177-191.
\item[{[23]}] Yurko V.A. Recovering variable order differential operators on star-type graphs
     from spectra. Differ. Uravneniya, 49, no.12 (2013), 1537-1548. (Russian); English transl. 
     in  Differ. Equations 49, no.12 (2013), 1490-1501.
\item[{[24]}] V.A. Yurko. Inverse problems on star-type graphs: differential operators of
     different orders on different edges. Central European J. Math. 12, no.3 (2014), 483-499.
\item[{[25]}] B.M. Levitan and I.S. Sargsyan. Introduction to Spectral Theory. AMS Transl. of Math.
     Monogr. 39, Providence, 1975.
\end{enumerate}

\begin{tabular}{ll}
Name:             &   Yurko, Vjacheslav  \\
Place of work:    &   Department of Mathematics, Saratov State University \\
{}                &   Astrakhanskaya 83, Saratov 410012, Russia \\
Present Position: &   Professor, Head of the Faculty of Mathematical Physics \\
Telephone:        &   (8452) 275526           \\
E-mail:           &   yurkova@info.sgu.ru
\end{tabular}

\end{document}